\documentclass[onecolumn,11pt]{article}
\usepackage[top=1.25in, bottom=1.25in, left=1.25in, right=1.25in]{geometry}
\usepackage{amsmath,amsfonts,amscd,amssymb}
\usepackage{graphicx}
\usepackage{overpic}
\usepackage{cancel}
\usepackage{rotating}
\usepackage{url}
\usepackage{caption}
\usepackage{color}
\usepackage{rotating}
\usepackage{multirow}
\usepackage{wrapfig}
\usepackage{mathtools}
\usepackage{subeqnarray}
\usepackage{setspace}
\usepackage{palatino} 
\graphicspath{{figures/}}

\usepackage[bottom,flushmargin,hang,multiple]{footmisc}
\usepackage{lipsum}
\newcommand\blfootnote[1]{%
  \begingroup
  \renewcommand\thefootnote{}\footnote{#1}%
  \addtocounter{footnote}{-1}%
  \endgroup
}

\newcommand{\x}{\hat{\mathbf{x}}}
\newcommand{\X}{\hat{\mathbf{X}}}
\newcommand{\kk}{\hat{\mathbf{k}}}
\newcommand{\K}{\hat{\mathbf{K}}}

\newcommand*\R{\mathbb{R}}

\title{\LARGE{\vspace{-.65in}\textbf{Smoothing and parameter estimation by soft-adherence to governing equations}}\vspace{-.25in}}
\author{\normalsize{Samuel H. Rudy$^{1*}$, Steven L. Brunton$^2$, J. Nathan Kutz$^1$}\\
\footnotesize{$^1$ Department of Applied Mathematics, University of Washington, Seattle, WA 98195, United States}\\
\footnotesize{$^2$ Department of Mechanical Engineering, University of Washington, Seattle, WA 98195, United States\vspace{-.8in}}
}
\date{}
\begin{document}
\maketitle

\blfootnote{$^*$ Corresponding author (shrudy@uw.edu).\\ \noindent \textbf{Python code:}  https://github.com/snagcliffs/SmoothingParamEst}
\vspace{-.25in}
\begin{abstract}
The analysis of high-dimensional dynamical systems generally requires the integration of simulation data with experimental measurements.  Experimental data often has substantial amounts of measurement noise that compromises the ability to produce accurate dimensionality reduction, parameter estimation, reduced order models, and/or balanced models for control.   Data assimilation attempts to overcome the deleterious effects of noise by producing a set of algorithms for state estimation from noisy and possibly incomplete measurements.  Indeed, methods such as Kalman filtering and smoothing are vital tools for scientists in fields ranging from electronics to weather forecasting.  In this work we develop a novel framework for smoothing data based on known or partially known nonlinear governing equations. The method yields superior results to current techniques when applied to problems with known deterministic dynamics.  By exploiting the numerical time-stepping constraints of the deterministic system, an optimization formulation can readily extract the noise from the nonlinear dynamics in a principled manner.  The superior performance is due in part to the fact that it optimizes global state estimates.  We demonstrate the efficiency and efficacy of the method on a number of canonical examples, thus demonstrating its viability for the wide range of potential applications stated above. \\

\vspace{-0.05in}
\end{abstract}



\section{Introduction}

The modern analysis of high-dimensional dynamical systems via data assimilation~\cite{evensen2009data} typically leverages a combination of simulation and experimental data, which is often obfuscated by noisy measurement~\cite{law2012evaluating,lee2015multiscale}.   Noise is well known to critically undermine one's ability to produce accurate low-dimensional diagnostic features, such as POD (proper orthogonal decomposition)~\cite{HLBR_turb} or DMD (dynamic mode decomposition)~\cite{kutz2016dynamic} modes, produce accurate forecasts and parameter estimations, generate reduced order models~\cite{benner2015survey}, or compute balanced truncation models for control~\cite{moore1981principal,willcox2002balanced,rowley2005model}.
The ubiquity of data imbued with noise has led to significant research efforts for filtering and parameter estimation in high dimensional chaotic dynamical systems, especially in application areas such as climate modeling, material science, fluid dynamics, and atmospheric sciences.
Past method for filtering include sequential Bayesian methods like the Kalman and particle filters~\cite{chui2017kalman}, and 3D and 4D-var \cite{courtier1994strategy}.  Parameter estimation for noisy data may be carried out by EM (expectation maximization) algorithms~\cite{moon1996expectation} alternating between state estimation and parameter estimation.  Central to each of these methods is balancing fidelity of state estimates to the prescribed dynamics and to measurements.
We present an alternative to these state-of-the-art mathematical techniques by exploiting the deterministic nature of the known underlying governing equations.  Specifically, we formulate an optimization whereby we enforce adherence to a time-stepping constraint for the deterministic system, thus separating the noise from the dynamics.

The filtering of time-series data is a well-developed field with a host of simple (e.g. bandpass filters) to sophisticated (e.g. the family of Kalman filters) mathematical architectures introduced over many decades.   For example, sequential Bayesian filters have gained immense popularity since their introduction by Kalman in 1960 \cite{Kalman1960jfe}.  Kalman filters obtain state estimates by tracking second order statistics of the error covariance for a trajectory with uncertain initial conditions, disturbances, and measurement error.  The Extended Kalman Filter (EKF) introduced methodology to filter signals from nonlinear dynamical systems \cite{maybeck1982stochastic} via linearizing the dynamics at each timestep to approximate changes in the error covariance.  However, the linearization used in the EKF may fail due to the first order expansion being inaccurate and may also be computationally challenging for large systems.  An improved Kalman filter for nonlinear dynamics called the Unscented Kalman Filter (UKF) was introduced in \cite{julier1997new} using the unscented transform to track error covariance by using the single timestep trajectories of a set of points chosen to mimic the initial state's error statistics.  The Ensemble Kalman Filter (EnKF) \cite{evensen1994sequential} was introduced for high dimensional dynamical systems and approximates statistics of the state estimate via a small ensemble of trajectories, generally numbering fewer than the dimension of the state space.  Each iteration of the Kalman filter uses Gaussian statistics to track the dynamics.  The particle filter was introduced for filtering problems on highly nonlinear systems where Gaussian statistics fail to accurately the error \cite{majda2010mathematical}.  Unfortunately, the particle filter exhibits poor performance in the high dimensional setting \cite{bengtsson2008curse}.  Recent work has developed particle filters for higher dimensional systems \cite{majda2014blended, chorin2004dimensional}, but full consideration of non-Gaussian statistics is considered in a low dimensional space.\\

State estimates in Kalman filters are conditioned on previous measurements and initial condition.  For fixed interval smoothing problems where full trajectory knowledge is available and state estimates do not need to be made in an online manner, Kalman smoother \cite{evensen2000ensemble} and the Rauch-Tung-Striebel (RTS) smoother introduce a backwards pass to the Kalman filter to condition state estimates on full trajectory measurements.  The RTS smoother is easily extended to nonlinear systems via the Unscented Rauch-Tung-Striebel filter (URTS) \cite{sarkka2008unscented} and Ensemble Rauch-Tung-Striebel filter (EnRTS) \cite{raanes2016ensemble}.  In contrast to sequential filters and smoothers, recent work \cite{rudy2018deep} has shown that large scale optimization problems may learn precise measurements of measurement error at every timestep, even when the dynamics are unknown.  In this work we explore the case where dynamics are known, possibly only up to a set of constants, and construct a method capable of both parameter estimation and highly accurate state estimation for non-linear and high dimensional systems, even where noise is non-Gaussian, correlated in time, or exhibits a constant offset.  In contrast to sequential Bayesian methods like the Kalman or particle filter, the proposed method optimizing a trajectory over the entire range of the data and enforces a soft-adherence to a known time-stepping scheme \cite{leveque2007finite}.  We demonstrate superior performance to the ensemble Kalman RTS smother on a selection of canonical problems.

\subsection{Contribution of this work}

This paper presents a novel method for smoothing experimental data where dynamics are known up to a set of coefficients.  We emphasize that the work does not share the generality of Bayesian sequential data assimilation methods such as the ensemble Kalman filter but does yield more accurate predictions on problems with deterministic dynamics.  Specifically, it optimizes over full time trajectories.  It is also possible to perform parameter estimation with the algorithm presented in this work without alternating between state and parameter estimation.  Finally, the method presented here is robust to non-Gaussian noise, stiff and highly nonlinear dynamics, and even noise with non-zero mean.  We believe this work forms a substantial contribution to any field where noisy measurements of smooth dynamics systems must be considered.

The paper is organized as follows:  In section \ref{methods} we formalize the problem statement, provide an overview of the mathematical justification for our methods, and introduce a novel computational method for smoothing data based on soft adherence to a time-stepping scheme.  In section \ref{results} we provide numerical results for a few standard test problems in data assimilation and draw comparisons to the Ensemble Rauch-Tung-Streibel Smoother.  In addition to Gaussian measurement error, error from Ornstein-Uhlenbeck processes, heavy tailed noise, and Gaussian error with non-zero mean are also considered.  Section \ref{discussion} contains a discussion of the method and remarks for further work.

\section{Methods}\label{methods}

This work leverages mathematical optimization to find state estimates for a measured time series by enforcing soft adherence to known dynamics.  We use known time-stepping architectures to construct measures of accuracy for estimated time series to our governing equations.  We also enforce  the state estimate to lie close to the measured data.  The proposed method is able to estimate the true state to high precision even in the case of temporally autocorrelated noise or noise with non-zero mean.  Section \ref{problem_formulation} provides a brief overview of the problems considered in this work and section \ref{comp_methods} introduces the proposed computational framework for denoising as well as a discussion of initialization and optimization.

\subsection{Problem formulation} \label{problem_formulation}

We consider observations ${\bf y}_j$ of a continuous process ${\bf x}$ given as,
\begin{equation}
\mathbf{y}_j = \mathbf{x}_j + \boldsymbol{\nu}_j,
\label{eq:noise_assumtion}
\end{equation}
where ${\bf x}_j \in \R^n$ denotes the value of ${\bf x}$ at time $t_j$ for $j = 1, \hdots , m$ and $\epsilon_j$ is measurement error.  We further presume that ${\bf x}$ obeys the autonomous differential equation,
\begin{equation}
\dot{\mathbf{x}} = f(\mathbf{x}),
\label{eq:dynamics}
\end{equation}
where $\mathbf{x}$ is the state and $f(\mathbf{x})$ the velocity.  Many modern techniques in data-assimilation and smoothing use sequential Bayesian filters that make forecast estimates of the state $\mathbf{x}_{j+1}$ from predictions of $\mathbf{x}_j$ and then assimilate that forecast using estimates of the error covariance and the measurement $\mathbf{y}_{j+1}$.  In contrast, this work seeks to estimate states $\mathbf{x}_j$ for each time step by balancing two criteria; 1) the estimate states should fit the prescribed dynamics and 2) the estimated state must be as close as possible to the measurements.

\subsection{Computational methods}\label{comp_methods}

The general form of a Runge-Kutta scheme~\cite{kutz2013data} is,
\begin{align}
\mathbf{x}_{j+1} &= \mathbf{x}_j + h_j \sum_{i=1}^s b_i \mathbf{k}_j^{(i)} \label{eq:rk_step} \\
\mathbf{k}_j^{(i)} &= f\left(\mathbf{x_j} + h_j\sum_{l=1}^s a_{il} \mathbf{k}_j^{(l)},\, t_j + c_i h_j\right) ,\label{eq:rk_intermediate}
\end{align}
where matrix $\mathbf{A} \in \R^{s\times s}$, and vectors $\mathbf{b}, \mathbf{c} \in \R^s$ define specific weights used in the timestepper.  Given $\mathbf{x}_j$ the solution to the system of equations in \eqref{eq:rk_step} \eqref{eq:rk_intermediate} gives the state $\mathbf{x}_{j+1}$ at the subsequent timestep and values of $f$ at intermediate stages.  If $a_{ij} = 0$ for $j\geq i$ then the method is explicit and there is no need to solve a set of algebraic equations, but these methods are generally less reliable for stiff problems.  For dense ${\bf A}$ it is possible to obtain $\mathcal{O}(h^{2s})$ accurate time stepping schemes.  For the remainder of this text, time dependence in $f$ will be omitted from notation, though including it will not change any analysis.\\

We note that for any pair of states $(\mathbf{x}_j ,\mathbf{x}_{j+1})$ from the ground truth solution for the denoising problem there exists a set of intermediate values of the velocity $\{\mathbf{k}_j^{(i)}\}_{i=1}^s$ such that equations \eqref{eq:rk_step} and \eqref{eq:rk_intermediate} are satisfied to within the order of accuracy for the timestepper being considered.  This observation forms the foundation for the methods considered in this work.  For any pair of state estimates $(\x_j, \x_{j+1})$ and intermediate values of the derivative $\{\kk_j^i\}_{i=1}^s$ we define cost functions that measure the ability of the intermediate values of the derivative to predict the next timestep, and the consistency of the intermediate values to their definition in \eqref{eq:rk_step} \eqref{eq:rk_intermediate}.  These are given by,
\begin{align}
L_1^{(j)}(\X, \K) &= \left\|\x_{j+1} - \left(\x_j + h_j \sum_{i=1}^s b_i \kk_j^{(i)}\right)\right\|_2^2 \label{eq:L1}\\
L_2^{(j,i)} (\X, \K) &= \left\| \kk_j^{(i)} -  f\left(\x_j + h_j\sum_{l=1}^s a_{il} \kk_j^{(l)}\right) \right\|_2^2  \label{eq:L2} ,
\end{align} 
where $\X = \{\x_j\}_{j=1}^m$ and $\K = \{(\kk_j^{(1)}, \hdots , \kk_j^{(s)})\}_{j=1}^m$.  For a time series of state estimates and estimates of intermediate values of the derivative, we define a global cost function evaluating the fidelity of the time series to the prescribed dynamics bu summing equations \eqref{eq:L1} and \eqref{eq:L2} over each timestep.  Since equations \eqref{eq:L1} and \eqref{eq:L2} do not reference the measurements $\mathbf{Y}$, we include a measure of closeness to the measured data as $g(\X-\mathbf{Y})$.  The cost function balancing accuracy to the timestepping scheme with closeness to data is given by,
\begin{equation}
L(\X, \K) = \sum_{j=1}^{m-1} \left( L_1^{(j)}(\X, \K) + \sum_{i=1}^s L_2^{(j,i)} (\X, \K) \right) + g(\X - \mathbf{Y}), \label{eq:cost_over_K} 
\end{equation}
where a suitable choice for $g(\cdot)$ may be the commonly used 2-norm $\lambda \| \cdot \|_2^2$ for Gaussian noise or the more robust 1-norm $\lambda \|\cdot \|_1$ for handling heavy tailed noise.\\

Optimization over \eqref{eq:cost_over_K} yields state estimates with considerable error.  Terms in \eqref{eq:cost_over_K} from the sum over \eqref{eq:L2} may be considerable different in magnitude than those from \eqref{eq:L2}.  Scaling the time variable by a suitable value of $\alpha$ so that $\tau= \alpha t$ and velocity is effectively scaled by $\alpha^{-1}$ provides a remedy for the imbalance but this is not practice in an application setting.  The appropriate scaling factor may be unclear since the magnitude of $f$ evaluated on the measurements may be significantly different than when evaluated on the optimal state estimate.  A more reliable fix is to optimize over intermediate values of the state rather than of the velocity.  For $f$ with non-singular Jacobian at each $\mathbf{x_j} + h_j\sum_{l=1}^s a_{il} \mathbf{k}_j^{(l)}$ the implicit function theorem tells us we can find an equivalent expression to \eqref{eq:rk_step},\eqref{eq:rk_intermediate} involving only values of $\mathbf{x}$.  Inverting $f$ at each $\mathbf{k}_j^{(i)}$ in eq. \eqref{eq:rk_intermediate} gives,
\begin{equation}
\mathbf{x}_j^{(i)} = f^{-1}(\mathbf{k}_j^{(i)}) = \mathbf{x_j} + h_j\sum_{l=1}^s a_{il} \mathbf{k}_j^{(l)} .\label{eq:x_j^i}
\end{equation}
Replacing each $\mathbf{k}_j^{(i)}$ in equations \eqref{eq:rk_step} and \eqref{eq:rk_intermediate} with $f(\mathbf{x}_j^{(i)})$ and applying $f^{-1}$ to both sides of \eqref{eq:rk_intermediate} we obtain a Runge-Kutta scheme where intermediate values are taken in the domain of $f$ rather than the range.
\begin{align}
\mathbf{x}_{j+1} &= \mathbf{x}_j + h_j \sum_{i=1}^s b_i f\left(\mathbf{x}_j^{(i)}\right) \label{eq:rk_step_x} \\
\mathbf{x}_j^{(i)} &= \mathbf{x_j} + h_j\sum_{l=1}^s a_{il} f\left(\mathbf{x}_j^{(l)}\right) .\label{eq:rk_intermediate_x}
\end{align}
Using equations \eqref{eq:rk_step_x} and \eqref{eq:rk_intermediate_x} we can follow the same process as before to construct a cost function over the set of variables $\X = \{\x_j\}_{j=1}^m$ and $\tilde{\mathbf{X}} = \{(\x_j^{(1)}, \hdots , \x_j^{(s)})\}_{j=1}^m$.  The resulting function for Runge-Kutta denoising is given by,
\begin{equation}
L_{RKD} (\X, \tilde{\mathbf{X}}) = \sum_{j=1}^{m-1} \left( \mathcal{L}_1^{(j)}\left(\X, \tilde{\mathbf{X}}\right) + \sum_{i=1}^s \mathcal{L}_2^{(j,i)}\left(\X, \tilde{\mathbf{X}}\right) \right) + g(\X - \mathbf{Y}), \label{eq:cost_over_X}
\end{equation}
where,
\begin{align}
\mathcal{L}_1^{(j)}\left(\X, \tilde{\mathbf{X}}\right) &= \left\|\x_{j+1} - \left(\x_j + h_j \sum_{i=1}^s b_i f\left(\x_j^{(i)} \right) \right)\right\|_2^2 \label{eq:L1_new}\\
\mathcal{L}_2^{(j,i)}\left(\X, \tilde{\mathbf{X}}\right) &= \left\| \x_j^{(i)} -  \left(\x_j + h_j\sum_{l=1}^s a_{il} f\left(\x_j^{(l)}\right)\right) \right\|_2^2  \label{eq:L2_new} .
\end{align} 
Note that \eqref{eq:cost_over_X} is very similar to \eqref{eq:cost_over_K} but where we are now optimizing over intermediate values of the state rather than the velocity.  It is also possible to include additional terms in \eqref{eq:cost_over_X} to penalize properties of the state estimate such as variation in derivatives to enforce smoothness.

Many state of the art algorithms for system identification from noisy data employ an expectation-maximization framework to alternate between smoothing and parameter estimation \cite{ghahramani1999learning}.  In contrast, the method considered in this work extends trivially to allow for learning state estimates and model parameters in a single minimization problem.  Where $f$ is known up to a set of parameters $\theta$ we replace $f$ by $\hat{f}_\theta$ in \eqref{eq:L1_new}-\eqref{eq:L2_new} and optimize over the state estimate with intermediate values and model parameters.

We minimize \eqref{eq:cost_over_X} using the quasi-Newton solver L-BFGS \cite{zhu1997algorithm}.  Derivatives of the cost function are evaluated using automatic-differentiation software implemented in the Tensorflow library in Python \cite{tensorflow2015}.  Initial estimates for the states are obtained by applying a naive smoother to the measured data and intermediate state estimates $\x_j^{(i)}$ are obtained by a linear interpolation between initial estimates of $\x_j$ and $\x_{j+1}$.   

\section{Numerical Examples}\label{results}

In this section we present several numerical examples of increasing complexity and dimension.  We test the denoising algorithm on the Lorenz 63, Lorenz 96, Kuramoto-Sivashinsky, and nonlinear Schr\"odinger equations.  We also provide examples of the denoising algorithm in cases where model parameters are unknown.  In each case the method is tested with white noise added to the state, but more complicated cases including noise drawn from an Ornstein-Uhlenbeck process, heavy-tailed noise, and severely biased noise are also considered.  We use the notation $\nu \sim N$ and let $\Sigma_N^2$ denote the measurement error covariance which we will report as a percentage of $\Sigma_X^2$, the variance of the state $X$.

\subsection{Lorenz 63} \label{L63}

The Lorenz 63 system given by \eqref{eq:L63} was first derived from a Galerkin projection of Rayleigh-Bernard convection and has since become a canonical teaching example for nonlinear dynamical systems \cite{Lorenz1963jas}.
\begin{equation}
\begin{aligned}
\dot{x} &= \sigma (y-x)\\
\dot{y} &= x(\rho - z) - y\\
\dot{z} &= xy - \beta z.
\end{aligned} \label{eq:L63}
\end{equation}

\begin{figure}
\centering
\includegraphics[width = \textwidth]{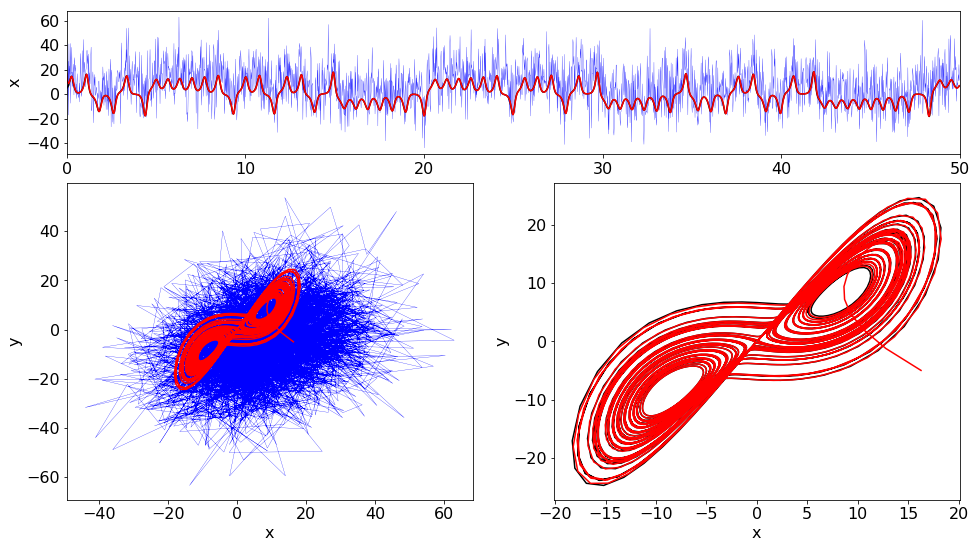}
\caption{Denoising results for Lorenz-63 system corrupted by white noise with non-zero mean.  $\nu \sim \mathcal{N}(\mu, \Sigma_X^2)$ with $\mu = (5,-5,-5)$.  RMSE = $0.397$}
\label{fig:Lorenz63_offset}
\end{figure}

We generated a test dataset by simulating the Lorenz 63 system for 2500 timesteps using timestep length $h = 0.02$, parameters $\sigma = 10$, $\rho = 28$, and $\beta = 8/3$, and initial condition $(5,5,25)$.  Several types of measurement error are added to the numerical solution to test the accuracy of the proposed method.  Figure \ref{fig:Lorenz63_offset} shows the data, state, and state estimate for the Lorenz 63 system for measurement noise distributed according to $\nu \sim \mathcal{N}(\mu, \Sigma_X^2)$ with $\mu = (5,-5,-5)$.  An initial transient exhibits considerable error but the remainder of the state estimate tracks the chaotic trajectory precisely.

\begin{figure}
\centering
\includegraphics[width = \textwidth]{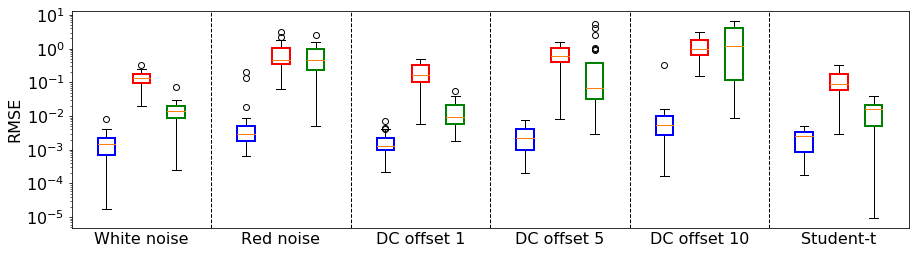}
\caption{Denoising results for Lorenz-63 system with $\sigma_\nu^2 = \sigma_X^2$ according to various distributions.  Blue: Runge-Kutta denoising with known dynamics, Red: Runge-Kutta denoising with unknown ($\mu$, $\sigma$, $\rho$), Green: Ensemble RTS Smoother with $N_e = 500$.}
\label{fig:Lorenz63_stats}
\end{figure}

A summary of results using root mean square error (RMSE) for the Lorenz 63 system with various noise distributions in shown in figure \ref{fig:Lorenz63_stats}.  In each case, noise is set to have variance equal to the data $X$.  We test the method for measurement error distributed according to mean zero Gaussian (white) noise, and Gaussian noise with non-zero mean having magnitude equal to 1, 5, or 10.  We also test the case where noise is autocorrelated in time (red noise) according to,
\begin{equation}
\begin{aligned}
&\nu_{j+1} = \rho \nu_j + \sqrt{1-\rho^2} \epsilon_j\\
&\epsilon_j, \nu_0 \sim \mathcal{N}(0, \Sigma^2)
\end{aligned}
\end{equation}
where $\rho = 0.75$.  The last case corresponds to an Ornstein-Uhlenbeck process and results in deviations from the true trajectory that are not remedied by moving average smoothing techniques.

In each case we test the method using the method introduced in this work with and without knowing the parameters $\rho$, $\sigma$, and $\beta$.  For comparison, we also obtain state estimates using the EnRTS smoother with initial state mean error covariance given by $(x_0, \Sigma_X^2)$, and process noise given by a mean zero normal distribution with covariance $dt^8I$ given the use of a fourth order timestepper.  We note that the inclusion of some error in from the timestepping scheme did not have significant effects on results.  The EnRTS smoother did not perform parameter estimation.  In each case the method presented in this work outperformed the EnRTS smoother when dynamics were known.  In the cases where the parameters were unknown, the proposed method resulted in a median RMSE comparable to the EnRTS smoother but generally having lower variance.

\subsection{Lorenz 96} \label{L96}

The Lorenz 96 system was introduced by Lorenz as a test model for studying predictability in atmospheric models \cite{lorenz1996predictability} and has since been shown to exhibit chaotic behavior \cite{karimi2010extensive}.  it is given by,
\begin{equation}
\dot{x}_i = (x_{i+1} - x_{i-2})x_{i-1} - x_i + F,
\label{eq:L96}
\end{equation}
where $F$ is a forcing parameter and $i = 1, \hdots, 40$ with period boundary conditions.  The degree of chaotic behavior is determined by $F$, with $F=8$ resulting in highly chaotic and $F=16$ resulting in turbulent behavior \cite{sapsis2013statistically}.

\begin{figure}
\centering
\includegraphics[width = \textwidth]{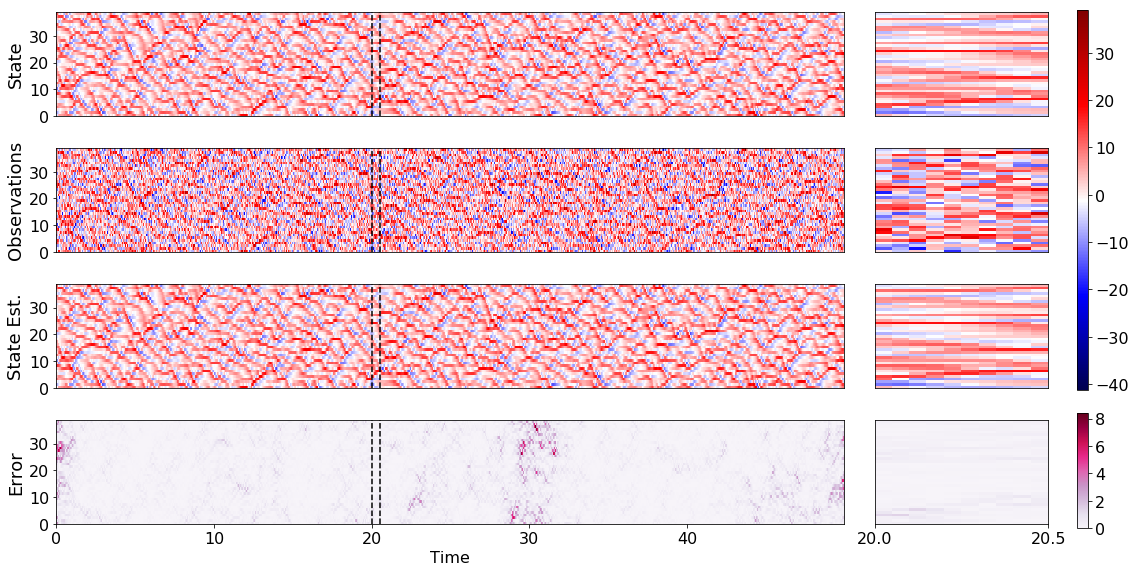}
\caption{Denoising and parameter estimation results for Lorenz-96 system with 100\% white noise.  $\hat{F} = 15.94$}
\label{fig:L96_param_est}
\end{figure}

\begin{figure}
\centering
\includegraphics[width = \textwidth]{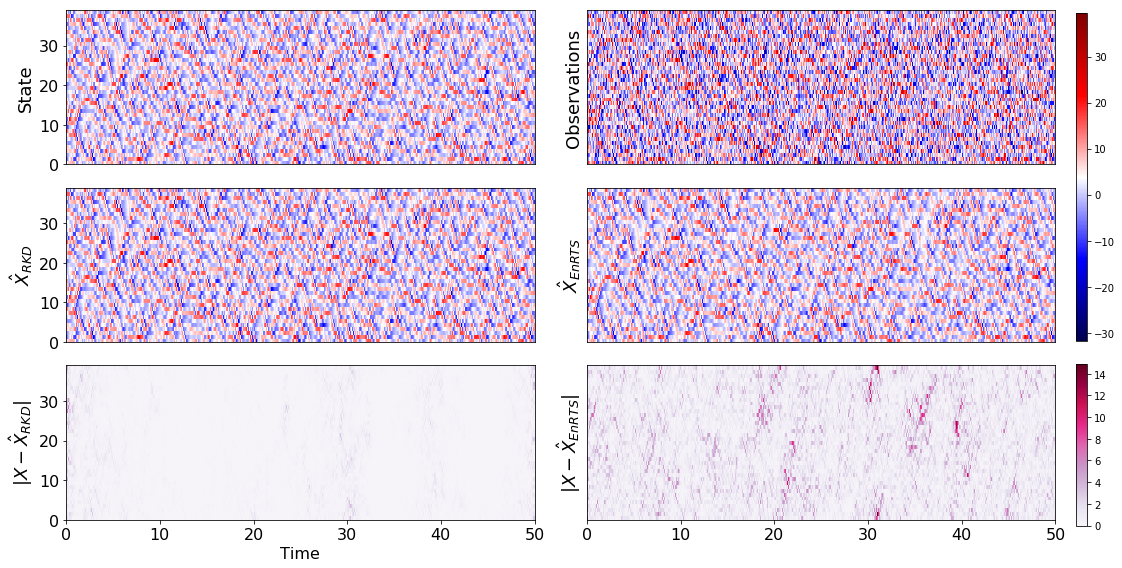}
\caption{Comparison between RKD and EnRTS for Lorenz-96 system with 100\% white noise.  Top left: $X$.  Top right: $Y$.  Center left: $X_{RKD}$. Center right: $X_{EnRTS}$. Bottom left: $|X-X_{RKD}|$. Bottom right: $|X - X_{EnRTS}|$.}
\label{fig:L96_surface_comparrison}
\end{figure}

Figure \ref{fig:L96_param_est} shows the results of the Runge-Kutta denoising algorithm on the Lorenz 96 system in the turbulent regime where the algorithm is not given the value of $F$.  A more detailed view of state, observations, and state estimate is shown for time $t \in [20,20.5]$.  The learned value of $F$ is $15.94$, $0.375$\% below the true value of $16$.

\begin{figure}
\centering
\includegraphics[width = \textwidth]{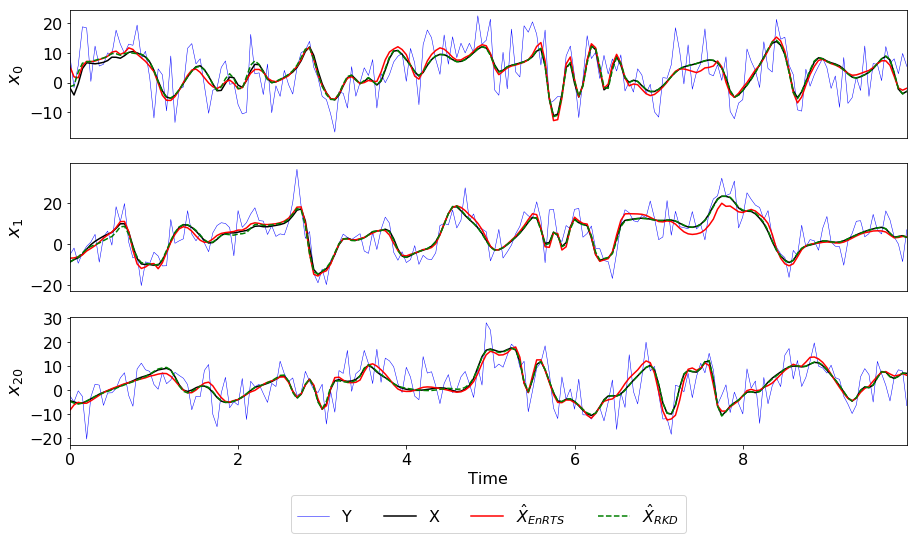}
\caption{Example of data and smoothed time series using RKD and EnRTS for Lorenz-96 system with 100\% white noise.}
\label{fig:L96_time_series_comparrison}
\end{figure}

\begin{figure}
\centering
\includegraphics[width = 0.667\textwidth]{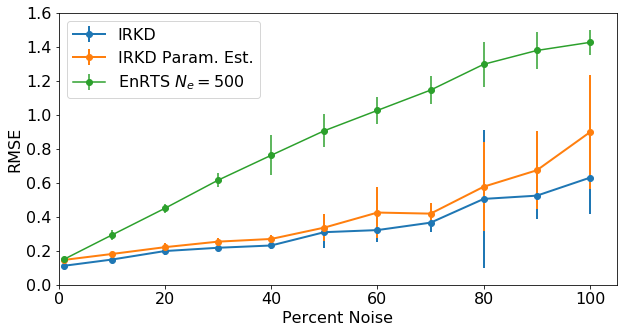}
\caption{Error for denoising Lorenz-96 system with 100\% white noise using EnRTS, RKD with known parameters, and RKD with parameter estimation.}
\label{fig:L96_comparrison_stats}
\end{figure}

Figures \ref{fig:L96_surface_comparrison} and \ref{fig:L96_time_series_comparrison} compare the results of the proposed algorithm to the EnRTS smoother on a single example of the Lorenz 96 system with known $F=16$.  Absolute error for the Runge-Kutta based smoothing is considerably lower than the EnRTS smoother.  The proposed method is also more accurate when it is required to learn the forcing parameter.  Figure \ref{fig:L96_comparrison_stats} shows the mean and standard deviation of the RMSE for Runge-Kutta based smoothing with and without known $F$, and the EnRTS smoother applied to the Lorenz 96 system in the turbulent regime across noise levels from $\sigma_N^2 = 0.01 \sigma_X^2$ to $\sigma_N^2 = \sigma_X^2$.  In each case the Runge-Kutta smoothing outperforms the EnRTS smoother even if $F$ is unknown.

Results for the Lorenz 96 system are only presented with measurement noise coming from a mean-zero Gaussian distribution.  If instead red noise was added to the state then both the proposed Runge-Kutta based smoothing and the EnRTS smoother failed to obtain accurate state estimates.

\subsection{Kuramoto Sivashinsky} \label{KS}

The Kuramoto Sivashinsky (KS) equation is a fourth order partial differential equation given by \eqref{eq:KS}.  It is considered a canonical example of spatiotemporal chaos  in a one-dimensional PDE \cite{hyman1986kuramoto, aceves1986chaos} and is therefore commonly used as a test problem for data-driven algorithms.  The KS equation is a particularly challenging case for filtering algorithms due to it's combination of high dimensionality and nonlinearity \cite{jardak2010comparison}.

\begin{equation}
u_t+uu_x+u_{xx}+u_{xxxx} = 0 \label{eq:KS} .
\end{equation}

\begin{figure}
\centering
\includegraphics[width=\textwidth]{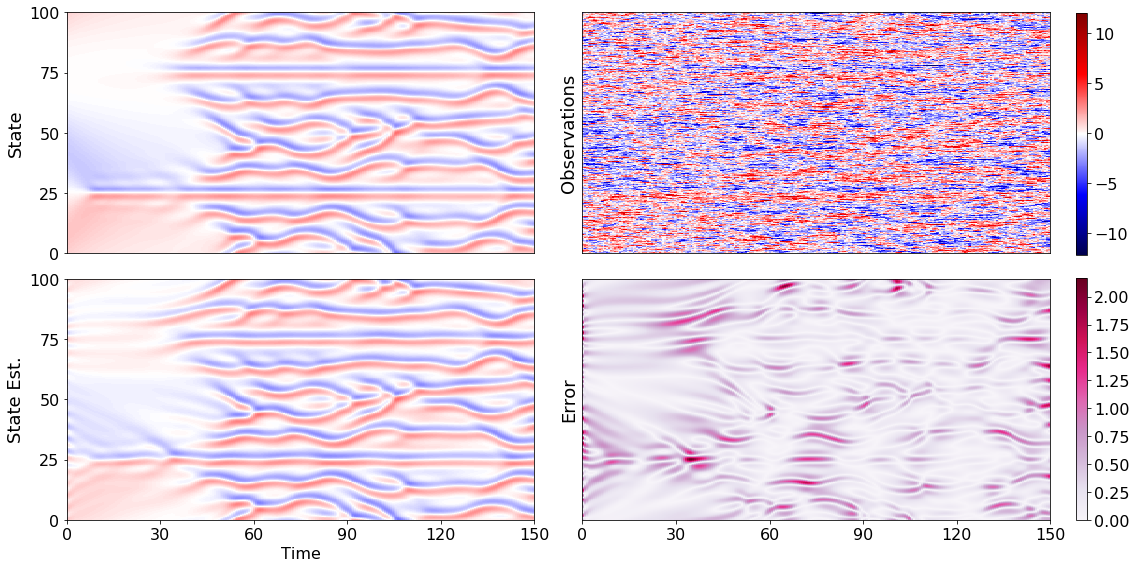}
\vspace{-.3in}
\caption{Results on Kuramoto Sivashinsky equation.  top left: $X$. Top right: $Y$ with $\sigma^2_N = 5\sigma^2_X$. Bottom left: $\hat{X}_{RKD}$. Bottom right: $|X - \hat{X}_{RKD}|$.}
\label{fig:ks_surface}
\end{figure}

\begin{figure}
\centering
\includegraphics[width = \textwidth]{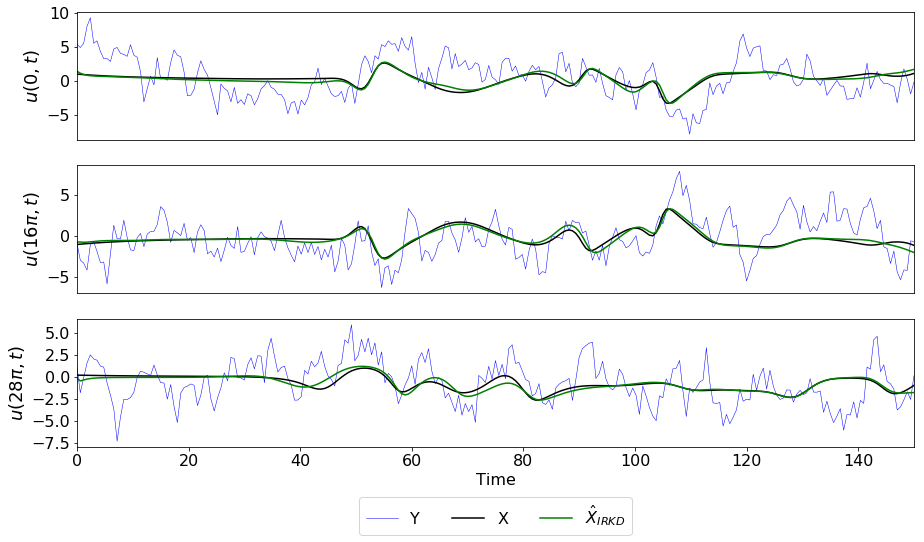}
\caption{Example of data and smoothed time series using RKD for KS system with 500\% red noise ($\rho = 0.75$).}
\label{fig:ks_series}
\end{figure}

We solve \eqref{eq:KS} on $[0, 32\pi]$ with periodic boundary conditions using the method outlined in \cite{kassam2005fourth} from $t=0$ to $t = 150$.  Figures \ref{fig:ks_surface} and \ref{fig:ks_series} show the results of the proposed method on the Kuramoto-Sivashinsky equation with red noise having $\Sigma_N^2 = 5 \Sigma_X^2$ and $\rho = 0.75$.

\subsection{Nonlinear Schr\"odinger Equation}

As an example of model reduction, we can consider the nonlinear Schr\"odinger (NLS) equation, which is a complex valued nonlinear PDE used as a model for nonlinear optics~\cite{kutz2006mode} and Bose-Einstein condensates~\cite{bronski2001bose}
\begin{equation}
iu_t + \frac{1}{2}u_{xx} + |u|^2u = 0 .
\label{eq:NLS}
\end{equation}
Breather solutions of the NLS equation, which are generated with initial conditions $u(x,0)=N\mbox{sech}(x)$ (here we consider $N=2$), exhibit a low-rank structure with 99\% of the energy contained in two modes making the equation a trivial example for the application of reduced order methods.  However, after corruption by small magnitude white noise, the rank necessary to capture 99\% of the energy in considerably higher.  We consider the effective rank of the breather solution to the NLS equation given a noisy dataset and after the proposed smoothing algorithm has been applied.  Table \ref{tab:nls_rank} summarizes the rank for a truncation preserving 99\% of the energy in the solution for noisy and smoothed datasets.  For clean data the effective rank is 2.  The algorithm proposed here for denoising can bring the rank of the system close to the ideal rank-2 evolution, thus providing a significantly improved reduced order model when compared to with the data which has not been denoised.

\begin{table}
  \begin{center}
    \caption{Rank for SVD truncation to 99\% energy of noisy NLS data.}
    \label{tab:nls_rank}
    \begin{tabular}{|l|c|c|c|c|c|c|}
      \hline
      Data & $1$ & $5$ & $10$ & $25$ & $50$ & $100$ \\
      \hline
      $\mathbf{Y}$ & 18 & 51 & 64 & 78 & 86 & 91 \\
      $\hat{\mathbf{X}}$ & 2 & 3 & 4 & 11 & 18 & 35 \\
      \hline
    \end{tabular}
  \end{center}
\end{table}

\begin{figure}
\centering
\includegraphics[width = \textwidth]{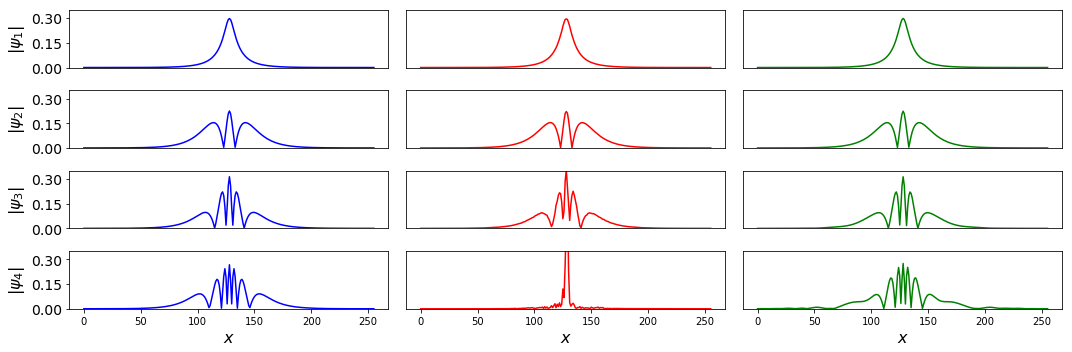}
\caption{Leading 4 principle vectors of NLS dataset for clean (blue), noisy (red), and smoothed (green) datasets for 100\% noise.}
\label{fig:NLS_odes}
\end{figure}

In the ideal case without noise, two modes clearly dominate the
behavior of the system.  These two POD modes are the first two
columns of the matrix ${\bf U}$ and are now used to approximate
the dynamics observed from full scale simulation.  In theory,
the two mode expansion takes the form~\cite{kutz2013data}
\begin{equation}
  u(x,t)=a_1(t) \phi_1(x) + a_2(t) \phi_2(x)
\end{equation}
where the $\phi_1$ and $\phi_2$ are the first two POD modes.  Inserting
this approximation into the governing NLS equation gives the reduced order model
\begin{subeqnarray}
\label{eq:pod2by2}
  && i {a_1}_t + \alpha_{11} a_1 + \alpha_{12} a_2 
  + \left( \beta_{111} |a_1|^2 + 2\beta_{211} |a_2|^2 \right) a_1
  \\
 && \,\,\,\,\,\,\,\,\,\,  + \left( \beta_{121} |a_1|^2 + 2\beta_{221} |a_2|^2 \right) a_2 + \sigma_{121} a_1^2 a_2^*
  + \sigma_{211} a_2^2 a_1^* = 0 \nonumber \\
   && i {a_2}_t + \alpha_{21} a_1 + \alpha_{22} a_2 
  + \left( \beta_{112} |a_1|^2 + 2\beta_{212} |a_2|^2 \right) a_1
  \\
 && \,\,\,\,\,\,\,\,\,\,  + \left( \beta_{122} |a_1|^2 + 2\beta_{222} |a_2|^2 \right) a_2
  + \sigma_{122} a_1^2 a_2^*
  + \sigma_{212} a_2^2 a_1^* = 0  \nonumber
\end{subeqnarray}
where $ \alpha_{jk} = ( {\phi_j}_{xx},\phi_k)/2$, $ \beta_{jkl} = ( |\phi_j|^2 \phi_k,\phi_l ) $ and $\sigma_{jkl}=( \phi_j^2 \phi_k^* ,\phi_l)$.  The initial values of the ROM are given by
$ a_1(0)= {(2\mbox{sech}(x),\phi_1)}/{(\phi_1,\phi_1)}$ and $ a_2 (0) = {(2\mbox{sech}(x),\phi_2)}/{(\phi_2,\phi_2)}$.  This gives a complete description of the two mode dynamics
predicted from the SVD analysis.

\begin{figure}
\centering
\includegraphics[width = \textwidth]{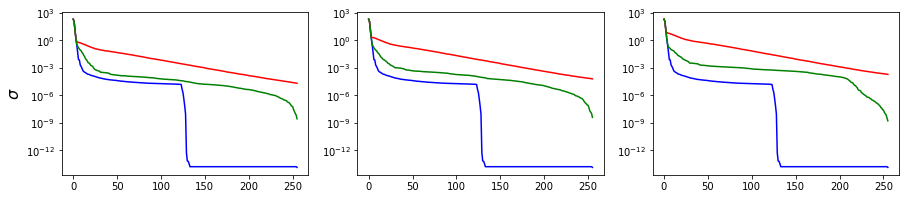}
\caption{Singular values of NLS data for clean (blue), noisy (red), and smoothed (green) data.  Left: 1\% noise, center: 10\% noise, right: 100\% noise.}
\label{fig:NLS_spectrum}
\end{figure}

The model (\ref{eq:pod2by2}) relies on accurate evaluations of the inner products between modes (both linear and nonlinear projections).  Figure~\ref{fig:NLS_odes} shows the first four modes of the system for clean data, the noisy data and the denoised data.  Although the first three modes of each is reasonable, the addition of noise has a significant influence on the 4th mode onwards.  A ROM built upon these modes performs poorly in characterizing the full fidelity model~(\ref{eq:NLS}).  Moreover, it suggests using a much higher rank ROM as shown in  Table~\ref{tab:nls_rank} and illustrated in Fig.~\ref{fig:NLS_spectrum}.  Denoising the data produces a ROM and rank truncation that are close to the ideal, noise-free data (at least up to 10\% noise).  This provides a much more accurate ROM and reconstruction of the spatiotemporal dynamics, enhancing forecasting possibilities with the ROM model.

\section{Discussion} \label{discussion}

In this work we have presented a novel technique for fixed interval smoothing and parameter estimation for high dimensional and highly nonlinear dynamical systems.  the proposed method forces adherence to known governing equations by minimizing the residual of a Runge-Kutta scheme while also pinning the state estimate to observed data.  While more limited than sequential Bayesian filters and smoothers, the method is robust to substantial measurement error and yields accurate results in several cases where Bayesian filters fail.  The proposed methodology is shown to be robust to non-Gaussian, offset, and temporally autocorrelated noise even in the high dimensional setting.

We believe that the optimization framework presented in this work suggests future research directions in the fields of system identification, parameter estimation, and filtering.  Directions could include extensions of the current work to include nonlinear or partial observation functions as well as considerations of stochastic dynamics.  The authors are committed to reproducible research and have made all code used in this work publicly available on GitHub at https://github.com/snagcliffs/SmoothingParamEst.

\section*{Acknowledgments} 
We acknowledge funding support from the Defense Advanced Research Projects Agency (DARPA- PA-18-01-FP-125). SLB further acknowledges generous funding from the Army Research Office (ARO W911NF-17-1-0306 and W911NF-17-1-0422) and JNK acknowledges support from the Air Force Office of Scientific Research (AFOSR FA9550-18-1-0200).  

\small
\small
\begin{spacing}{.5}
\bibliographystyle{plain}
\bibliography{denoising_bibliography}
\end{spacing}

\normalsize

\section*{Appendix A: Optimization over velocity space}\label{time_scaling}

In section \ref{comp_methods} we suggest that optimization over incremental values of the state using equation \eqref{eq:cost_over_X} is superior to optimizing over incremental velocities as is standard in Runge-Kutta methods.  One possible reason for the difference is the potential for large differences in magnitude between $\mathbf{x}$ and $\dot{\mathbf{x}}$.  If so, the scale difference could be remedied by scaling time by $\tau = \alpha t$ so that $\|\mathbf{x} \| \approx \|\mathbf{x}_\tau\| = \alpha^{-1} \|\dot{\mathbf{x}}\|$.\\

\begin{figure}[h]
\centering
\includegraphics[width=0.667\textwidth]{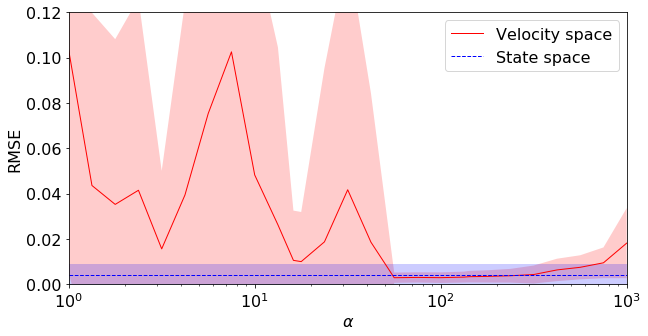}
\caption{Root mean square error for the solution of \eqref{eq:cost_over_K} across 25 trials for various values of $\alpha$ in red.  RMSE for \eqref{eq:cost_over_X} shown in blue. }
\label{fig:time_scaling_error}
\end{figure}

Figure \ref{fig:time_scaling_error} investigates the performance of the proposed smoothing technique using equation \eqref{eq:cost_over_K} (red) and equation \eqref{eq:cost_over_X} (blue).  Indeed, there exist scalings $\alpha$ such that optimization over velocity space performs comparatively with optimization over state space.  However, finding the optimal scaling may not be possible for a naive application of the smoothing technique where the true state is unknown.  Figure \ref{fig:time_scaling_error} indicates that for too small of a time dilation $\alpha$ the method's performance exhibits substantial variability across trials.

\section*{Appendix B: Measure of fidelity to measured data}\label{regularization}

In equation \eqref{eq:cost_over_X}, measurements are only included in the term $g(\hat{\mathbf{X}}- \mathbf{Y})$.  We use either the $L^1$ norm or squared $L^2$ norm as a means of pinning the state estimate to the measured data with a constant $\lambda$ dictating the importance of this difference relative to the terms measuring fidelity to the dynamics.  Without the term $g(\hat{\mathbf{X}}-\mathbf{Y})$, the data $\mathbf{Y}$ would only appear in the optimization problem as an initial guess for the solver.  A reasonable hueristic for setting $\lambda$ would be to balance the expected loss due to measurement noise $g(\hat{\mathbf{N}})$ for some reasonable guess of the noise $\mathbf{N} = [\boldsymbol{\nu}_1, \hdots, \boldsymbol{\nu}_m]$ with the expected loss due to timestepper error, which can be roughly approximated by $m(s+1)dt^p$ for an order-$p$ timestepping scheme having $s$-steps.  However, the method is generally very robust to changes in $\lambda$ so a careful analysis may not be necessary.\\

\begin{figure}
\centering
\includegraphics[width=0.667\textwidth]{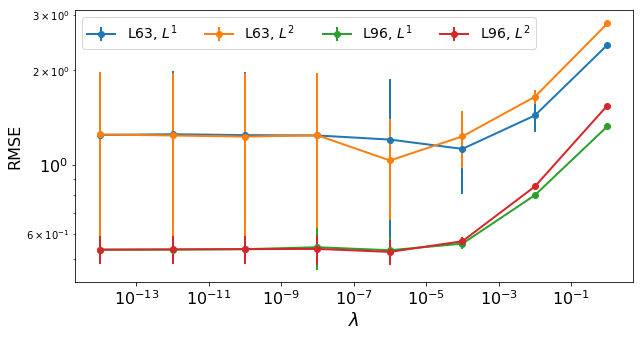}
\caption{Root mean square error for the solution of \eqref{eq:cost_over_X} using $g(\cdot) = \lambda \|\cdot \|_1$ or $g(\cdot) = \lambda \|\cdot \|_2^2$ for various $\lambda$ averaged over 10 random noise samples un the Lorenz 63 and Lorenz 96 systems.  $\lambda$ on the order of $10^{-8}-10^{-4}$ provides small advantage over $\lambda=0$ solution.}
\label{fig:regularizer_sensitivity}
\end{figure}

Figure \ref{fig:regularizer_sensitivity} shows the sensitivity of the algorithm to changes in the weighting $\lambda$ of $g(\cdot)$ in the cases where either an $L^1$ or $L^2$ norm is used for the Lorenz 63 and Lorenz 96 systems.  Surprisingly, the method generally performed well even when $\lambda = 10^{-18}$, essentially negating the role of data as anything other than an initial guess.  However, this result may not be general and slight improvements in performance were observed for $\lambda \in [10^{-8}, 10^{-4}]$.  Each example shown in the paper used $\lambda = 10^{-8}$.\\

\end{document}